\documentclass{article}
\usepackage{latexsym}
\usepackage{amsmath}
\usepackage{amssymb}

\newcommand{\ev}{{e\!v}}

\newcommand{\R}{{\mathbb R}}

\newcommand{\Q}{{\mathbb Q}}

\newcommand{\Z}{{\mathbb Z}}

\newcommand{\CP}{{\mathbb CP}}

\newcommand{\Hh}{{\mathcal H}}

\newcommand{\Ss}{{\mathcal S}}
\newcommand{\Ll}{{\mathcal L}}
\newcommand{\Uu}{{\mathcal U}}

\newcommand{\oMm}{{\overline{\mathcal M}}}

\newcommand{{\TS}}{{\Tilde {\rm Symp}}}
\newcommand{{\THH}}{{\Tilde {\rm Ham}}}
\newcommand{\p}{{\partial}}
\newcommand{\al}{{\alpha}}

\newcommand{\be}{{\beta}}

\newcommand{\Om}{{\Omega}}
\newcommand{\om}{{\omega}}
\newcommand{\eps}{{\varepsilon}}

\newcommand{\ga}{{\gamma}}
\newcommand{\Ga}{{\Gamma}}
\newcommand{\ka}{{\kappa}}
\newcommand{\la}{{\lambda}}
\newcommand{\La}{{\Lambda}}
\newcommand{\si}{{\sigma}}

\newcommand{\Symp}{{\rm Symp}}
\newcommand{\Diff}{{\rm Diff}}
\newcommand{\Ham}{{\rm Ham}}

\newcommand{\Flux}{{\rm Flux}}

\newcommand{\Si}{{\Sigma}}

\newcommand{\MS}{{\medskip}}
\newcommand{\SSK}{{\smallskip}}

\newcommand{\NI}{{\noindent}}
\newcommand{\proof}[1]{\noindent{\bf Proof#1:\  }}
\newcommand{\QED}{\hfill${\bf QED}$\medskip}

\newtheorem{theorem}{Theorem}[section]
\newtheorem{thm}[theorem]{Theorem}
\newtheorem{cor}[theorem]{Corollary}
\newtheorem{lemma}[theorem]{Lemma}

\newtheorem{prop}[theorem]{Proposition}

\begin{document}

\title{Cohomological properties of 
ruled symplectic structures}
\author{Fran\c{c}ois Lalonde\thanks{Partially supported by NSERC grant 
OGP 0092913
and FCAR grant ER-1199.} \\ Universit\'e du Qu\'ebec \`a Montr\'eal
\\ (flalonde@math.uqam.ca) \and Dusa McDuff\thanks{Partially
supported by NSF grants DMS 9704825 and 0072512.} \\ State University of New York
at Stony Brook \\ (dusa@math.sunysb.edu)
}

\date{October 12, 2000}
\maketitle

\MS

\section{Introduction}

Donaldson's work on Lefschetz pencils has shown that,
after a slight perturbation of the symplectic form and a finite
number of blow-ups, 
any closed symplectic manifold $(M,\om)$ can be expressed as a singular 
fibration  with generic fiber a smooth
codimension 2 symplectic submanifold. Thus fibrations play a fundamental
role in symplectic geometry. It is then natural to study 
smooth (non-singular)  ruled symplectic manifolds $(P, \Om)$, i.e.
symplectic manifolds where $P$ is the total space of a smooth fiber bundle 
$M \hookrightarrow P
\to B$ and the symplectic form $\Om$ is such that its restriction $\om$ 
to each $M$-\/fiber is nondegenerate. 

  In this survey, we present some of the recent results obtained by the
authors and Leonid Polterovich
in~\cite{LMP,McD,LM-Topology} that show that 
bundles endowed with such structures have interesting stability and
cohomological  properties. For instance, 
under certain topological conditions on the
base, the rational cohomology of $P$ necessarily splits as the tensor product
of the cohomology of $B$ with that of $M$. 

   As we describe below in \S~\ref{se:char},
the bundle $M \hookrightarrow P \to B$ corresponding to a ruled
symplectic manifold, has the group $\Ham(M,\om)$ of Hamiltonian
diffeomorphisms of $M$ for structure group if the base $B$ is simply
connected.
They can therefore be divided into two
classes: those whose structural group belongs to a finite dimensional
Lie subgroup of the group $\Ham(M)$ of Hamiltonian
diffeomorphisms of $M$, and those whose structural group
is genuinely infinite dimensional.  The first case belongs
to the realm of classical symplectic geometry and its study 
can be carried out by topological and Morse-theoretic methods,
whereas the second case belongs to symplectic topology and
requires analytic methods such Gromov-Witten
invariants and versions of quantum homology which are relevant
in the context of bundles. Our results tend to confirm that
many properties of the first case still hold in the second one. 
They therefore fit in with the principle mentioned by
Reznikov~\cite{Rez}
that the group of symplectomorphisms behaves like a
Lie group. \MS 

Our main conjecture is

\begin{quote}  The rational cohomology of a bundle
$M \hookrightarrow P \to B$, with $\Ham(M,\om)$ for structure
group, splits as the tensor product of the cohomology of 
the fiber with that of the base. In particular, the same splitting
occurs for any ruled symplectic  manifold $(P,\Om)$ over a simply 
connected base.
\end{quote}

  The section \S~\ref{se:finite} of this survey presents the results 
related to the finite dimensional case, that were already
known (Blanchard, Deligne, Kirwan, Atiyah-Bott). These results prove
the above conjecture in the finite dimensional case and
provide a good intuition of what should hold in 
general.  The next section,
\S~\ref{se:infinite}, presents our results in the infinite dimensional case.
These
are established under some restrictions on the topology of the base, and 
lead to corollaries that apply to the general
theory of ruled symplectic manifolds. 

It is still unclear
whether our conjecture holds in the full generality in which it is stated, 
that is to say when the structure group is infinite dimensional and when
no restrictions are placed on the topology of the base and fiber.
However, our methods
are sufficiently general to yield all known results
in the finite dimensional case.
More details about many of the results presented here, together with some 
applications to the action of $\Ham(M,\om)$ on $M$, can be found 
in~\cite{LM-Topology}.

\MS

Before getting into these questions, it is useful to give a
characterisation of ruled symplectic manifolds that holds
when the base $B$ is simply connected.

\section{Characterizing ruled symplectic manifolds}\label{se:char}

A symplectic manifold $(P,\Om)$ is said to be ruled if $P$ is a 
(locally trivial) fiber bundle over some base manifold $B$ and the 
restriction of $\Om$ to each fiber is nondegenerate.
 It turns out that there is a close relation between such
manifolds and {\it Hamiltonian bundles}.

A fiber bundle $M\to P\to B$ is said to be {\it symplectic}
(resp. {\it Hamiltonian}) if, for some symplectic form $\om$ on $M$,
 its structural group
reduces to the group of symplectomorphisms $\Symp(M, \om)$ of $(M,\om)$
(resp. to the group $\Ham(M, \om)$ of Hamiltonian diffeomorphisms of
$(M,\om)$).
In both cases, each fiber $M_b = \pi^{-1}(b)$ is equipped with a well defined
symplectic form $\om_b$ such that $(M_b,\om_b)$ is symplectomorphic to 
$(M, \om)$.    It is easy to see that every ruled  symplectic
manifold may be given the structure of a symplectic fiber bundle.
However, it is not so obvious that
when the base is simply connected this bundle can be taken 
Hamiltonian.

Here is a geometric criteria for a symplectic bundle
to be Hamiltonian, i.e. for the structural group to reduce to $\Ham(M, \om)$. 
Recall that the group  $\Ham(M, \om)$ is a normal subgroup of the identity
component $\Symp_0(M,\om)$ of the group of symplectomorphisms, and fits into the
exact sequence
$$
\{id\} \to \Ham(M, \om) \to \Symp_0(M, \om)\stackrel{\Flux}\to 
H^1(M,\R)/\Ga_{[\om]}  \to
0,
$$
where $\Ga_{[\om]}$ is a countable group called
 the flux group. In particular, $\Ham(M,\om)$ is
connected,
which means that a Hamiltonian bundle is trivial over the $1$-skeleton of the
base.    The
following proposition was proved in~\cite{MS}~Thm.~6.36 by a somewhat
analytic argument and in \cite{LM-Topology} by a more geometric
one.

\begin{prop}\label{prop:hamchar}  A symplectic bundle
  $\pi: P\to B$ is 
Hamiltonian
 if and only if the
following conditions hold:\begin{itemize}
 \item[(i)]  the restriction of $\pi$ to the $1$-skeleton $B_1$ of $B$
is  symplectically trivial, and
\item[(ii)] there is
a cohomology class $a\in H^2(P,\R)$ that restricts to $[\om_b]$ on $M_b$.
\end{itemize}
\end{prop}

When the spaces are smooth manifolds,  a construction due to 
Thurston shows that the existence of the extension class $a$ in 
$(ii)$ above is equivalent to the existence of a closed form 
$\tau$ on $P$ that  extends the family of forms $\om_{b}$ on the 
fibers.  Such a form $\tau$ gives rise to a connection on the 
bundle $P\to B$ whose horizontal distribution consists of the 
$\tau$-orthogonals to the fibers.  In this language, the previous 
result can be reformulated as follows.

\begin{prop}\label{prop:hamchar2}
A  symplectic
bundle  $\pi: P\to B$ is Hamiltonian if and only if the forms $\om_b$ 
on the fibers have  a
closed extension $\tau$ such that the  holonomy of the corresponding
connection
around any loop in $B$ 
lies in the identity component $\Symp_0(M)$ of $\Symp(M)$. 
\end{prop}

   When the base is itself a symplectic manifold, we can
add that form $\tau$ to a sufficiently large multiple of the pull-back
of the form on the base, which is then easily seen to be nondegenerate
in all directions. This proves:

\begin{cor} A Hamiltonian bundle over a
symplectic base is a ruled symplectic manifold.
Conversely, a ruled symplectic manifold over a simply connected
base is a Hamiltonian bundle.
\end{cor}

  However the two categories are not the same:  
the trivial Hamiltonian bundle $M\times 
S^{1}$ is not ruled, while the Kodaira--Thurston manifold $X$ that 
fibers over $T^{2}$ with nontrivial monodromy is ruled but not Hamiltonian.

\section{The finite dimensional case} \label{se:finite}

We now consider the results in the finite dimensional case
that were obtained in various contexts (Deligne in the
algebraic case, Kirwan for circle actions using localisation techniques
and Atiyah--Bott for torus actions). 

   Suppose that
the structural group of the Hamiltonian bundle
$(M,\om) \hookrightarrow P \to B$ can be reduced to
a compact Lie group $G \subset \Ham(M,\om)$. This means that there is
a representation of the group $G$ in the group of all
Hamiltonian diffeomorphisms of $M$, that is to say a Hamiltonian
action of $G$ on $M$. In this context, it is enough to consider
the universal Hamiltonian $G$-bundle with fiber $M$
$$
M \;\longrightarrow\; M_G = EG\times_G M \;\longrightarrow\; BG.
$$
The cohomology of $P = M_G$ is known as the equivariant cohomology
$H_G^*(M)$ of $M$. Kirwan showed in \cite{K} that if $G$ is the circle,
 the rational cohomology of the 
total space $M_G$ is isomorphic to the tensor product of the rational
cohomology of $M$ and the rational cohomology of $BG$ (we will refer
to this in the sequel by saying that the bundle $M_G$ is {\it c-split},
i.e. it splits cohomologically). This was proved by localisation techniques.
A different argument was given by Atiyah--Bott  in \cite{AB}, establishing 
the same result when $G$ is a torus of arbitrary dimension.
The result for a general 
compact Lie group $G$ follows by a more or less standard argument, 
starting either from
Atiyah--Bott's result or from our infinite dimensional results (see next section).
Here is the argument:

\begin{prop}\label{prop:funct2}
If $G$ is a compact connected Lie group that acts in a Hamiltonian way on $M$,
then any bundle $P\to B$ with fiber $M$ and structural group $G$ is c-split. 
In particular, 
$$
H^*_G(M) \cong H^*(M)\otimes H^*(BG).
$$  
\end{prop}

\proof{}  It is enough to prove the second statement since 
$$
M_G = EG\times _GM \longrightarrow  BG
$$
is the universal bundle.
Every compact Lie group $G$ is the image of a homomorphism
$T\times H \to G$, where the torus $T$ maps onto the identity component of the
center of $G$ and $H$ is the semi-simple Lie group corresponding to the
commutator subalgebra $[{\rm Lie} (G), {\rm Lie} (G)]$ in the Lie algebra
${\rm Lie}(G)$.  It is easy to see that this  homomorphism induces a
surjection on rational homology $BT \times BH \to BG$.  Therefore, we may
suppose that $G = T\times H$.  
Let $T_{max} = (S^1)^k$ be the
maximal torus of the semi-simple group $H$. Then 
the induced map on cohomology $H^*(BH) \to H^*(BT_{max}) =
\Q[a_1,\dots, a_k]$ takes $H^*(BH)$ bijectively onto the set of polynomials in
$H^*(BT_{max})$ that are invariant under the action of the Weyl group, by the
Borel-Hirzebruch theorem. Hence  the maps $BT_{max} \to BH$
and $BT \times BT_{max} \to BG$
induce a surjection on homology.   Therefore the desired result
follows from part (ii) of the following  lemma:

\begin{lemma}\label{le:funct} Consider a commutative diagram
$$
\begin{array}{ccc}
P' & \to & P\\
\downarrow & &\downarrow\\
B'& \to & B
\end{array}
$$
where $P'$ is the induced bundle.  Then:
\begin{itemize}
\item[(i)]  If $P\to B$ is c-split so is $P'\to B'$.\SSK

\item[(ii)] {\bf (Surjection Lemma)} If $P'\to B'$ is c-split and $H_*(B')\to
H_*(B)$ is surjective, then $P\to B$ is c-split.\end{itemize}
\end{lemma}
\proof{} $(i)$:  Use the fact that, by the Leray-Hirsch theorem, 
$P\to B$ is c-split if and only if the map $H_*(M)\to
H_*(P)$ is injective.

\NI
$(ii)$:  The induced map on the $E_2$-term of the cohomology spectral 
sequences is
injective.  Therefore the existence of a nonzero differential in 
the spectral sequence $P\to
B$ implies that the corresponding differential for the pullback bundle 
$P'\to B'$ does not vanish
either. \QED

   Smooth projective bundles constitute an interesting 
special case of the above 
proposition. The result in this case can be derived from the Deligne spectral
sequence, or more generally by the following argument due to Blanchard \cite{Bl}.  
Let's call a smooth fiber bundle $M \hookrightarrow P \to B$ {c-Hamiltonian} if
there is a class $a\in H^2(P)$ whose restriction $a_M$ to the
fiber  $M$ is c-symplectic, i.e. $(a_M)^n\ne 0$ where $2n = \dim(M)$.
 Recall that  a closed manifold $M$ is said to
satisfy the hard Lefschetz condition with respect to the class
$a_M\in H^2(M, \R)$ if the maps
$$
\cup (a_M)^k: H^{n-k}(M,\R) \to H^{n+k}(M,\R),\quad 1\le k\le n,
$$ 
are isomorphisms.   In this case, elements in $H^{n-k}(M)$ that vanish when
cupped with $(a_M)^{k+1}$ are called primitive, and the cohomology of
$M$ has an additive basis consisting of elements of the form $b\cup (a_M)^\ell$
where $b$ is primitive and $\ell \ge 0$.

\begin{prop}{\bf (Blanchard~\cite{Bl})}   Let  $M\to P\to B$ be a c-Hamiltonian
bundle  such that $\pi_1(B)$ acts trivially on $H^*(M, \R)$.  If in addition
$M$ satisfies the hard Lefschetz condition with respect to the c-symplectic
class $a_M$, then the bundle c-splits. \end{prop}
\proof{}
The proof is by contradiction.  Consider the Leray spectral sequence in 
cohomology  and suppose that
$d_p$ is the first non zero differential.  Then, $p\ge 2$ and
 the $E_p$ term in the spectral sequence is isomorphic to
the $E_2$ term and so can be identified with the tensor product $H^*(B)\otimes
H^*(M)$. Because of
the product structure on the spectral sequence, one of the differentials $d_p^{0,i}
$ must be nonzero.  So there is $b\in E_p^{0,i} \cong H^i(M)$ such that
$d_p^{0,i}(b) \ne 0$.  We may assume that $b$ is primitive (since these
elements together with $a_M$ generate $H^*(M)$.)
Then $b \cup a_M^{n-i} \ne 0$ but $b \cup a_M^{n-i+1} = 0$.
 
We can write $d_p(b) = \sum_j e_j\otimes f_j$ where $e_j \in H^*(B)$
and $f_j \in H^\ell(M)$ where $\ell < i$.  Hence $f_j \cup a_M^{n-i+1}\ne 0$
for all $j$ by the Lefschetz property.  Moreover, because the $E_p$ term is a
tensor product
 $$
(d_p(b))\cup a_M^{n-i+1} = \sum_j e_j \otimes (f_j \cup a_M^{n-i+1}) \ne 0.
$$
But this  is impossible since this element is the image via $d_p$ of the 
trivial element $b\cup a_M^{n-i+1}$.

Here is another, perhaps easier, argument.
Suppose $d = d_p$ is the first nonvanishing differential.
It vanishes on $H^i(M)$ for 
$i < p$ for reasons of dimension.    Therefore, by the Lefschetz property
it must vanish on 
$H^{2n-i}(M)$ for these $i$.  But then it has to vanish on $H^i(M)$ for
$p \le i < 2p$.  For if not, take $b$ in such $H^i(M)$ such that $d(b) \neq 0$.
By Poincar\'e duality there is $c \in H^{2n-j}(M)$ for $0\le j < p$ such that $d(b)
\cup c \ne 0$.  But $b\cup c = 0$ for reasons of dimension, and so $b\cup
d(c) \ne 0$, a contradiction.  It follows that  $d$ vanishes on
$H^{2n-j}$ for $p\le j < 2p$.  Now consider the next block of $i$:
$2p \le i < 3p$ and so on.
\QED

\MS

Another fundamental question about Hamiltonian bundles is that of their
stability under small perturbations of the symplectic form on the
fiber.  If the bundle $M\to P\to B$ has structure group $\Ham(M, \om)$ and $\om'$ 
is some nearby form, an elemenatry argument (given in 
Lemma~\ref{le:stab0} below) shows that it can be given the 
structure group $\Symp_{0}(M, \om')$.  However, it is not at all obvious 
whether the latter group can be reduced to $\Ham(M, \om')$.  If
this reduction is possible for all $\om'$ close to $\om$, 
 the original Hamiltonian bundle is said to be {\it stable}.
When the stuctural group is a compact Lie group, this
clearly boils down to the following statement.

\begin{thm} \label{thm:stability} $\,$ {\rm \bf (Hamiltonian 
stability.)} Let
$(M, \om)$ be a closed symplectic manifold, and let $\iota: G \to
\Ham(M, \om)$ be a continuous homomorphism defined on a compact Lie group.
Then, for each perturbation $\om'$ in some sufficiently small
neighbourhood $\Uu$ of $\om$ in the space of all symplectic forms on $M$,
there is a continuous homomorphism  
$$
\iota': G \to
\Ham(M, \om')
$$
that varies continuously as the form $\om'$ varies 
in $\Uu$.
\end{thm}
\proof{}
 We begin
with a well-known averaging argument.
Define $\tau'$ to be the average of the forms $\iota(g)^*(\om')$, i.e.
set $$
\tau'(v,w) = \int_G \iota(g)^*(\om')(v,w)\, d\mu_G, \qquad v,w\in T^*(M),
$$
where $d\mu_G$ is Haar measure.  Since $G$ is compact and
$\iota(g)^*(\om) = \om$ for all $g\in G$, $\tau'$ is a symplectic form 
whenever $\om'$ is sufficiently close to $\om$.  Moreover, it is easy to
see that $\iota(g)^*(\tau') = \tau'$ for all $g\in G$.  Thus $\iota$ maps $G$ into
$\Symp(M, \tau')$.  But, since $\iota(G)$ is also contained in the connected
group $\Ham(M,\om)$, the elements of $G$ must act trivially on $H^2(M)$.
Therefore $\tau'$ is cohomologous to
$\om'$ and hence equals $f^*(\om')$ where $f$ is the time $1$ map of
some isotopy $f_t$ (again assuming that $\om'$ is
 sufficiently close to $\om$.) 
Thus, the homomorphism
$$
\iota':\quad g\mapsto f\iota(g) f^{-1}
$$
takes $G$ to $\Symp_0(M,\om')$.  

It remains to show that this homomorphism  $\iota'$ has image in $\Ham(M, 
\om')$,  in other words that the composite homomorphism
$$
G\stackrel{\iota'}{\longrightarrow} \Symp_0(M,\om')
\stackrel{\Flux} {\longrightarrow} H^{1}(M, \R)/\Ga_{[\om']}
$$ 
is zero.  Since the target is an abelian group
and $G$ is compact, it suffices to consider the case when $G$ is the circle. 
Observe also that
 $\iota' =\iota_{\tau'}$ and $\iota
 = \iota_{\om}$ are the same when considered as maps into $\Diff(M)$.
 Therefore the same
 vector field $X$ generates both $S^1$-actions and we just need to show
that the  closed $1$-form
 $\al_{\tau'}$ which is $\tau'$-dual to $X$
 is actually exact. But each of
 its Morse-Bott singular sets have the same topology and index as does the
corresponding singular set for the  $\om$-dual $\al_{\om}$ of
 $X$,  since the underlying $X$ is the same and $\tau'$ is
arbitrarily close to  $\om$. Since $\al_\om$ is exact, all these indices
are even dimensional. Now if $x_0$ is a global
 minimum of a primitive $H_{\om}$ of $\al_{\om}$, then it is a local minimum
 of any local primitive of $\al_{\tau'}$. If $\ga$
 is any loop based at $x_0$, it is easy to see that
 $\int_{\ga} \al_{\tau'}$ must vanish, because otherwise the map $t\mapsto
 \int_{\ga(0)}^{\ga(t)} \al_{\tau'}$ would take different values at $t=1$
and $t=0$ and then a standard minimax argument over the loops
homotopic to  $\ga$ would yield a loop $\ga_0$ such that the maximum value
of the function
 $$
[0,1] \to \R,\quad t\mapsto \int_{\ga_0(0)}^{\ga_0(t)} \al_{\tau'}
$$
 would  occur at a critical point
 of $\al_{\tau'}$ of odd index, a contradiction. \QED

\section{The general case:  quantum homology and geometric induction on
ruled symplectic manifolds} \label{se:infinite}

   The methods developed in the finite dimensional case do not apply
when the structural group $\Ham(M,\om)$ of a Hamiltonian bundle 
$M \hookrightarrow P \to B$ does not retract to a finite dimensional
subgroup. Note that, even in examples as simple as $S^2 \times S^2$
endowed with a generic symplectic structure,
the Hamiltonian group does not retract to any Lie subgroup (see \cite{AMcD}). 

  It turns out that quantum homology and the Gromov-Witten invariants
can be used to derive
properties concerning the ordinary homology of a Hamiltonian bundle in full
generality, without any hypothesis on the structure group or on the 
symplectic manifold $(M,\om)$. It leads to a proof of the stability
of ruled symplectic structures and other corollaries
that we explain in the next sections.

  Basically,  pseudoholomorphic techniques (quantum homology,
GW-invar\-iants) apply when the fiber of the Hamiltonian bundle
is any compact symplectic manifold
and the base is
a 2-sphere. This means that the bundle is given by a loop in the
group $\Ham(M)$, which is generally not autonomous (i.e. it is 
a continuous map $S^1 \to \Ham(M)$ that need not be a homomorphism).
In this case, as we explained in \S~\ref{se:char}, there is a natural symplectic 
structure $\Om$ on the total space and we can first equip that space
with an almost complex structure compatible both with $\Om$ and with 
its restriction
to fibers (this can be done in general when the base is any symplectic manifold)
and then study the pseudo-holomorphic sections of
$$
\pi: (P, \Om, J) \to (S^2,j)
$$
where $j$ is the complex structure of $S^2$.  The moduli space
of these sections can be used to pair the quantum homology of
the $M$\/-fiber at say the north pole with the one at the south pole.
The main point is that this pairing is nondegenerate, i.e. it induces
an isomorphism between the two quantum homologies. It is then easy
to show that this implies that the GW-invariant attached to any 
ordinary cycle in the fiber at the north pole, {\it considered as a cycle in the
total space} $P$,  cannot vanish if it does not vanish as a cycle in the
fiber. In other words, the rational homology of the fiber injects inside
the homology of $P$, which implies by the Leray-Hirsch theorem
that $P$ is c-split. 

    We can then extend  these results to Hamiltonian
bundles defined over more general bases than $S^2$ by 
using three types of arguments:
\MS

\NI
1) Analytic arguments
that constitute a non-trivial generalisation of the analytic methods 
used over $S^2$. Basically, these arguments prove the c-splitting of Hamiltonian
fiber bundles over any base $B$ that has enough $J$-holomorphic rational curves;
more precisely, one asks that that $B$ has at least one non-vanishing rational 
GW-invariant in some class $A \in H_2(B)$
of the form $n(pt, pt, c_1, \ldots, c_k;A)$ where $k \ge 0$ and the $c_i$'s
are any homology classes in $B$ (see Proposition~\ref{prop:split}.).\SSK

\NI
(2) Geometric arguments needed to iterate bundles; they lead to a proof
that a Hamiltonian bundle over some base $P$ c-splits if $P$ is itself a 
Hamiltonian bundle $M \hookrightarrow P \to B$ over a simply connected base $B$
and if all Hamiltonian bundles having $M$ or $B$ as base c-split. 
\SSK

\NI
(3) Topological arguments based on properties of spectral sequences 
of symplectic bundles.\MS

\subsection{Bundles over $S^2$}

   We begin by explaining this in more detail in the case of a Hamiltonian
bundle over $S^2$ (this was proved in the semi-monotone case
by the authors in collaboration
with Polterovich in~\cite{LMP} and in the general case in~\cite{McD}.)  

   There is a correspondence between loops in the group of symplectic
diffeomorphisms and  
symplectic bundles over $S^2$ with fiber
$(M, \om)$.   The
correspondence is  given by assigning 
to each symplectic loop $\phi_{t \in
[0,1]}$ in $\Symp_0(M)$ the bundle $(M, \om) \to P_{\phi} \to S^2$
obtained by gluing a copy of $D_2^+ \times M$ 
with another $D_2^- \times M$
along their boundary in the following way:  
$$
(2\pi t, x) \mapsto  (-2\pi t, \phi_t(x)).
$$
(Here $D_2$ is the closed disc of radius $1$ of the plane.)
In what follows we
always assume that the base $S^2$ is oriented, and with orientation induced
from $D_2^+$. Note that this correspondence can be reversed: given a
symplectic bundle over the oriented $2$-sphere together with an
identification of one fiber with $M$, one can reconstruct the homotopy class
of $\phi$.

An important topological tool for the study of such bundles 
is the {\it Wang exact sequence}: 
$$
 ...\to H_{j-1}(M,\Z) \stackrel{{\partial_{\phi}}}{\to} H_j(M,\Z)
\stackrel{i} \to H_j(P_\phi,\Z) \stackrel{\cap [M]} {\to} H_{j-2}(M,\Z) \to ...
$$
  This sequence can be easily
derived
from the exact sequence of the pair $(P_\phi, M)$, where  
 $M$ is identified with a fiber of $P_\phi$.  The important point for us is
that the boundary map $H_{j-1}(M)\to H_j(M)$ is precisely the trace
homomorphism $\p_\phi$ that assigns to each cycle $a$
the cycle $\{\phi_t(a) \, | \, t \in [0,1] \}$.  Thus $\p_\phi$ vanishes
exactly when the inclusion $i$ is injective or, equivalently, when the
restriction map $\cap [M]$ is surjective. By the Leray-Hirsch theorem,
this is equivalent to the c-splitting of the bundle.

\begin{thm} \label{thm:sphere-c-splits}  Let $\phi$ be a Hamiltonian loop
on a closed symplectic manifold $(M,\om)$.
Then the homomorphism $i:H_*(M,\Q) \to H_*(P_{\phi},\Q)$
is injective; that is to say, the bundle c-splits.
\end{thm}

 We now briefly explain how the proof of this theorem proceeds.
Since $p:P_{\phi} \to S^2$ is a Hamiltonian
bundle it carries a natural deformation class  of symplectic forms
given by the weak coupling construction. 
Recall that {\it the coupling class} $u_\phi \in
H^2(P_{\phi},\R)$ is the (unique) class whose top power vanishes, and whose
restriction to a fiber coincides with the cohomology class of the fiberwise
symplectic structure. Let $\tau$ be a positive generator
of $H^2(S^2,\Z)$. The deformation class above consists of
symplectic forms $\Om$ which represent
the cohomology
class of the form $ u_\phi + \ka\, p^*\tau$  ($\ka >>0$)
and extend the fiberwise symplectic structure.   It is always possible to
choose $\Om$ so that it is a product with respect to the given product
structure near the fibers $M_0$ at $0\in D_2^+$ and $M_\infty$ at 
$0\in D_2^-$.

Besides the coupling class $u_\phi$,
the total space $P_{\phi}$ carries another canonical
second cohomology class 
$$
c_\phi = c_1(TP_\phi^{{\rm vert}})  \in H^2(P_{\phi},\R)
$$
that is defined to be the first Chern
class of the vertical tangent bundle. 

 Both  classes $u_\phi, c_\phi$ behave well under
compositions of loops.  More precisely, consider two elements $\phi, \psi\in
\pi_1(\Ham(M,\om))$ and their composite $\psi*\phi$.  This can be
represented either by the product $\psi_t\circ\phi_t$ or by the concatenation
of loops.   It is not hard to check that the bundle $P_{\psi*\phi}$ can be
realised as the fiber sum $P_\psi\#P_\phi$ obtained as follows.
Let $M_{\phi,\infty}$ denote the fiber at $0\in D_2^-$ in $P_\phi$ and
$M_{\psi,0}$ the fiber at $0\in D_2^+$ in $P_\psi$. Cut out open product
neighborhoods of each of  these fibers and then glue the
complements by an orientation reversing symplectomorphism of the
boundary.   The resulting space may be realised as
$$
D_2^+\times M\;\cup_{\al_{\phi, -1}}\;S^1\times [-1,1]\;\cup_{\al_{\psi,1}} 
\;D_2^-\times M,
 $$
where
$$
{\al_{\phi, -1}}(2\pi t,x) = (2\pi t, -1,\phi_t(x)),\quad 
{\al_{\psi, 1}}(2\pi t, 1,\psi_t(x)) = (2\pi t,x),
$$
and this may clearly 
be identified with  $P_{\psi*\phi}$. Set
$$
V_\phi = D_2^+\times M\;\cup\;S^1\times [-1,1/2), \quad 
V_\psi = \;S^1\times (-1/2,1]\;\cup 
\;D_2^-\times M.
$$

The next lemma follows imediately from  the construction of the coupling form
via symplectic connections.

\begin{lemma} The classes $u_{\psi*\phi}$ and $c_{\psi*\phi}$ are
compatible with the decomposition $
P_{\psi*\phi}= V_\psi\cup V_\phi$ in the sense that their
 restrictions to $V_\psi\cap V_\phi = (-1/2,1/2)\times S^1 \times M$ equal
the pullbacks of $[\om]$ and $c_1(TM)$.
\end{lemma}

  We now explain the proof of Theorem~\ref{thm:sphere-c-splits} (see
\cite{LMP,McD} for more details).

\bigskip

\NI
{ \large  Seidel's maps $\Psi_{\phi,\si}$}
\MS

We start with the definition of the quantum cohomology ring of $M$.
To simplify our formulas we will
denote the first Chern class $c_1(TM)$ of $M$ by $c$.

Let $\La$ be the usual (rational) Novikov ring of the group $\Hh =
H_2^S(M,\Z)/\!\!\sim$  with valuation $\om(.)$ where $B\sim B'$
if $\om(B-B') = c(B-B') = 0$, and   let $\La_R$ be
the analogous (real) Novikov ring based on the group $\Hh_R =
H_2^S(M,\R)/\sim$. Thus the elements of $\La$ have the form $$
\sum_{B\in \Hh} \la_B e^B
$$
where for each $\ka$ there are only finitely many nonzero
$\la_B\in \Q$ with $\om(B) < \ka$,
 and the elements of
$\La_R$ are
$$
\sum_{B\in \Hh_R} \la_B e^B,
$$
where 
$\la_B\in \R$ and there is a similar finiteness condition.\footnote
{
In~\cite{Sei} Seidel works with a simplified version of the Novikov
ring $\La$ where the coefficients $\la_B$ affecting $e^B,\; B \in \Hh$,
are elements of $\Z/2\Z$. However, his results extend
in a staightforward way to the case of rational coefficients
by taking into account orientations on the moduli spaces
of pseudo-holomorphic curves. Let us emphasize that in our
definition of $\La_R$ not only the coefficients $\la_B$ are real,
but also the exponents $B$ belong to a real vector space $\Hh_R$.}
Set 
$QH_*(M)
= H_*(M)\otimes\La$ and $QH_*(M,\La_R)
= H_*(M)\otimes\La_R$.  Then $QH_*(M)$ 
is $\Z$-graded with $\deg(a\otimes e^B) = \deg(a) - 2c(B)$.  
It is best to
think of
$QH_*(M,\La_R)$ as $\Z/2\Z$-graded with 
$$
QH_{\ev} =  
H_{\ev}(M)\otimes\La_R, \quad QH_{{\rm odd}} =  
H_{{\rm odd}}(M)\otimes\La_R.
$$
 With respect to the quantum intersection product 
both versions of quantum homology are graded-commutative
rings with unit $[M]$. Further, the units in $QH_{\ev}(M,\La_R)$ form a group 
$QH_{\ev}(M,\La_R)^\times$ that acts on   $QH_*(M, \La_R)$ by
quantum multiplication. 

Now we describe how Seidel arrives at an
action of the loop $\phi$ on the quantum homology of $M$.
Denote by $\Ll$ the space of
contractible loops in the manifold $M$.  
Fix a constant loop $x_* \in \Ll$, and
define a covering $\Tilde \Ll$ of $\Ll$ with the base point
$x_*$ as follows. Note first that a path between $x_*$ and
a given loop $x$ can be considered as a $2$-disc $u$ in $M$
bounded by $x$. Then the covering $\tilde \Ll$ is defined by saying that
two paths are equivalent if the $2$-sphere $S$ obtained by gluing
the corresponding discs has $\om(S) = c(S) = 0$. 
Thus the covering group
of $\Tilde \Ll$ coincides with the abelian group $\Hh$.

Let $\phi = \{\phi_t\}$ be a loop of Hamiltonian diffeomorphisms.
Because the orbits $\phi_t(x), t\in [0,1],$ of $\phi$ are contractible
(see~\cite{LMP1}), one can define a mapping $T_{\phi}: \Ll \to \Ll$ 
which takes
the loop $\{x(t)\}$ to a new loop $\{\phi_t (x(t))\}$. Let
$\Tilde T_{\phi}$
be a lift of $T_{\phi}$ to $\tilde \Ll$. To choose
such a lift one should specify a homotopy class of paths
in $\Ll$ between the constant loop and an orbit of $\{\phi_t\}$.
It is not hard to see that in the language 
of symplectic bundles this choice of lift
corresponds to a choice of an equivalence class $\si$ of 
sections of $P_{\phi}$,
where two sections are identified if their values under
$c_\phi$ and $u_\phi$ are equal.
Thus the lift can be labelled $\Tilde T_{\phi,\si}$.

Recall now that the Floer homology $HF_*(M)$  can be
considered as the Novikov homology of the action functional
on $\Tilde \Ll$.
Therefore $\Tilde T_{\phi, \si}$ defines a natural automorphism $(\Tilde
T_{\phi,\si})_*$ of $HF_*(M)$. Further, 
if $\Phi: HF_*(M) \to HQ_*(M)$ is the
canonical isomorphism constructed in 
Piunikhin--Salamon--Schwartz~\cite{PSS},
there is a corresponding 
automorphism $\Psi_{\phi,\si}$ of $QH_*(M)$ given by
$$
\Psi_{\phi,\si}\; = \;\Phi \circ ({\Tilde T_{\phi,\si}})_* \circ \Phi^{-1}.
$$
This gives rise to an action of the group of all pairs $(\phi,\si)$ on $QH_*(M)$.

Seidel shows in ~\cite{Sei} that when $M$ satisfies a suitable semi-positivity condition the map
$$
\Psi_{\phi,\si}: QH_*(M)\to QH_*(M)
$$
is in fact induced by quantum multipication
by an element of $QH_{\ev}(M)^\times$ that is formed from the moduli space of all
$J$-holomorphic sections of $P_\phi$.  
 In our work we went backwards. 
 We gave a  new
definition of the maps $\Psi_{\phi,\si}$ that does not  mention Floer
homology, and  proved that they are isomorphisms by a direct gluing
argument. 
Besides being easier to work with for general $M$,
our version of the definition
 no longer restricts us to using the 
coefficients $\La$ via the covering $\Tilde
\Ll\to \Ll$.  Instead we will consider the extension $\La_R$, which
 allowed us to  define an action of the group $\pi_1(\Ham)$ itself
(see \cite{LMP}).

\MS

 Let $\Om$ be  a
symplectic form  on $P_\phi$ that extends $\om$ and is in the
natural deformation class $u_\phi + \ka\,p^*(\tau)$.  As above, define an
equivalence relation on the set of homology classes of sections
of $P_{\phi}$ 
by identifying two such classes if their values under $c_\phi$ and $u_\phi$
are equal.   
Then, given a loop of Hamiltonian diffeomorphisms $\phi$ on $M$,
and an equivalence class of sections $\si$ of $P_{\phi}$
with $d = 2c_\phi(\si)$, define a $\La$-linear map
$$
\Psi_{\phi,\si}: \;\; QH_*(M) \to QH_{*+d}(M)
$$
as follows.   For $a \in
H_*(M,\Z)$, $\Psi_{\phi,\si}(a)$ is  the class in $QH_{*+d}(M)$ whose
intersection with $b\in H_*(M,\Z)$ is given by: 
$$
  \Psi_{\phi,\si}(a) \cdot_M
b = \sum_{B\in \mathcal H} n(i(a),i(b); \si + i(B))e^B.
$$
Here, we have written $i$ for the homomorphism
$H_*(M) \to H_*(P)$ and $\cdot_M$ for the linear extension to $QH_{\ast}(M)$
of the usual  intersection
pairing on $H_*(M,\Q)$.  Thus
$a\cdot_M b = 0$ unless $\dim(a) + \dim (b) = 2n$ in which case it is the
algebraic number of intersection points of the cycles.
Further, $n(v,w;D)$ denotes the Gromov--Witten invariant
which counts isolated $J$-holomorphic stable curves 
in $P_\phi$ of genus $0$
and two marked points that
represent the equivalence class $D$ and whose marked points go
through  given generic
representatives of the classes $v$ and $w$ in $H_*(P_{\phi},\Z)$. 
More precisely, one  defines $n(v, w; D)$ to be
the intersection of the virtual moduli cycle
$$
\ev: \oMm_{0,2}\,\!^\nu(P_\phi, J, D) \to P_\phi\times P_\phi,
$$
that consists of all perturbed 
$J$-holomorphic genus $0$ stable maps that lie in
class $D$ and have $2$ marked points, with a generic representative of
the class
$v \otimes w$ 
in $P_\phi\times P_\phi$. This definition is well understood
provided 
$M$  is spherically monotone or has minimal
spherical  Chern number~\footnote
{
The minimal spherical Chern number $N$
is the smallest nonnegative 
integer such that the image of $c = c_1(TM)$ on
$H_2^S(M)$ is contained in $N\Z$.  The weakly exact
case $N = 0$ is also tractable by these standard methods.
}
 $N \ge n-1$. In the general case, one uses a version 
of the virtual moduli cycle for Gromov-Witten invariants that is adapted 
to the fibered structure: see~\cite{McD}.

  Note finally that,
by Gromov compactness, there are for each
given energy level $\ka$  only finitely many homology classes $D$
with $\om(D-\si) \le \ka$ that are represented by $J$-holomorphic curves
in $P_\phi$.  Thus
 $\Psi_{\phi,\si}(a)$ satisfies the finiteness condition for
elements  of $QH_*(M,\La)$.

Since $n(i(a), i(b); D) = 0$ unless $
2c_\phi(D) + \dim (a) + \dim (b) = 2n,
$
we have 
$$
\Psi_{\phi,\si}(a) = \sum a_B\otimes e^B,
$$
where 
$$
\dim (a_B) = \dim(a) + 2c_\phi(D)  = \dim (a) + 2c_\phi(\si) + 2c(B).
$$
 Observe also that
$$
\Psi_{\phi,\si + B} = \Psi_{\phi,\si}\otimes e^{-B}.
$$

When $M$  is spherically monotone or has minimal
spherical  Chern number at least $n-1$
the following two results are proved by Seidel~\cite{Sei}.  The general
case is established in \cite{McD}.  

\begin{lemma} \label{le:constant} If $\phi$ is the constant loop 
$*$ and $\si_0$ is the flat
section $pt\times S^2$ of $P_* = M\times S^2$ then $\Psi_{*,\si_0}$ is the
identity map. 
\end{lemma}

\begin{prop} \label{prop:gluing} Given  sections $\si$ of $P_\phi$ and  $\si'$ of
$P_\psi$ let $\si'\#\si$ be the union of these sections in the fiber sum
$P_\psi\#P_\phi = P_{\psi*\phi}$.  Then 
$$
\Psi_{\psi,\si'}\circ \Psi_{\phi,\si} = \Psi_{\psi*\phi,
\si'\#\si}.
$$
\end{prop}

The main step in the proof of these statements
is to show that when
calculating the Gromov-Witten invariant $n(i(a),i(b);D)$
via the intersection between the virtual moduli cycle
and the class $i(a) \otimes i(b)$
we can assume the following:

\NI
--- the representative of $i(a) \otimes i(b)$ has the form
$\al\times \be$ where $\al, \be$ are cycles 
lying in the fibers of $P_\phi$;
 
\NI
--- the intersection occurs with elements in the top stratum of
$\oMm_{0,2}\,\!^\nu(P_\phi, J, D)$ consisting of sections of $P_\phi$.

\NI
In the semi-monotone case, Lemma~\ref{le:constant} is 
then almost immediate
\footnote{The proof of the first lemma is surprisingly hard in the general 
case. 
 The difficulty lies in showing that invariants in classes $A + B$
with $B\ne 0 \in H_2(M)$ do not 
contribute.  The reason is that such curves can never be isolated:
 they are graphs, and reparametrizations of the map to $M$ give 
rise to families of graphs. However, to see this in the general case
 involves 
constructing a virtual moduli cycle that is invariant under an $ S^{1}$-action.
See \cite{McD}.}, and  
Proposition~\ref{prop:gluing} can be proved by the
well-known  gluing techniques of~\cite{RT} or~\cite{MS1}.

\begin{cor}  $\Psi_{\phi,\si}$ is an isomorphism for all loops
$\phi$ and sections $\si$. 
\end{cor}

With this in hand, we can establish
the proof of Theorem~\ref{thm:sphere-c-splits} in the following
way. The Gromov-Witten invariants are linear in each variable.
Thus if $i(a)=0 $ for some $a \neq 0$, then $\Psi_{\phi,\si}(a) = 0$,
a contradiction with the fact that $\Psi_{\phi,\si}$ is an 
isomorphism. \QED
\MS

\subsection{General stability} \label{ss:General-stability}

   It turns out that the fact that all Hamiltonian bundles over $S^2$
are c-split is enough to establish the stability of general Hamiltonian 
bundles over any base. This is what we now explain.

Let  $\pi:P\to B$ be a symplectic bundle with closed fiber $(M,\om)$
and compact base $B$.  
 Moser's homotopy
argument implies that this bundle has the following stability property. 

\begin{lemma}\label{le:stab0} There is an open neighborhood 
$\Uu$ of $\om$ in the space $\Ss(M)$ of all symplectic forms
such that $\pi:P\to B$ may be naturally considered as a symplectic bundle 
with fiber $(M,\om')$ for all $\om'\in U$.  
\end{lemma}
\proof{} First recall that for every symplectic structure $\om'$ on $M$ there
is a Serre fibration 
$$
\Symp(M,\om') \longrightarrow \Diff(M) \longrightarrow  \Ss_{\om'},
$$
where $\Ss_{\om'}$ is the space of all symplectic structures on $M$ that are
diffeomorphic to $\om'$.\footnote
{Note that one uses Moser's argument to prove
that this is a Serre fibration.
}
At the level of classifying spaces, this gives a homotopy fibration
$$
\Ss_{\om'} \hookrightarrow  B \Symp(M, \om') \longrightarrow B \Diff(M).
$$
Any fiber bundle $P\to B$ with fiber $M$ is
classified by a map $B \to B\,\Diff(M)$, and isomorphism classes of symplectic
structures on it with fiber $(M, \om')$ correspond to homotopy classes of
sections of the associated bundle $W(\om') \to B$ with fiber $\Ss_{\om'}$.
If we are given a finite set of local
trivializations $T_i: \pi^{-1}(V_i) \to V_i\times M$, 
the transition functions of $P$ are maps 
$\phi_{ij}: V_i\cap V_j \to \Diff(M)$.
We can of course choose them with values in $\Symp(M,\om)$ if the $T_i$'s
are chosen compatible with the $\om$\/-structure on the bundle, but this 
is not necessary. Then the {\it
same} transition functions can be used to define the bundle 
$\Ss_{\om'} \hookrightarrow W(\om') \to B$, whatever the
symplectic structure $\om'$ may be.
 
Therefore, we are given a section $\si$ of $W(\om)$ and our task is to show 
how to construct from it a smooth family of sections $\si_{\om'}$ of the bundles
$W(\om')$ for all $\om'$ near $\om$.
Let $\si_i$  be the restriction of $\si$ over $V_i$.  Then $(T_i)_* \si_i$ is a
smooth map $V_i\to \Ss_\om$ (constant and equal to $\om$ if the $\phi_{ij}$ are
chosen in $\Symp(M, \om)$).  For each
$\om'$ near $\om$ and
$b\in V_i$ consider the symplectic form 
$$
\si_i' (b) =(T_{i})_*^{-1}((T_{i})_* (\si_i(b)) + \om' - \om)
$$
on $M_b$ (this is simply $(T_{i})_*^{-1}( \om')$ if the transition
functions have values in $\Symp(M,\om)$).  This is
cohomologous to $\om'$, and because of the openness of the symplectic nondegeneracy condition
will belong to the fiber of $W(\om')$ at
$b$, whenever $\om'$ is sufficiently close to $\om$.  Therefore, $\si_i'$ is
a section of $W(\om')$ over $V_i$. 
Morever, if 
$b\in \cap_{i\in I} V_i$ for some index set $I = I_b$, the convex hull of the
forms $\si_i'(b),
i\in I_b,$ will consist entirely of symplectic forms isotopic to $\om'$ and so 
will lie in the fiber of $W(\om')$ at $b$, again provided that 
$\om'$ is sufficiently close to $\om$.  Hence if
$\rho_i$ is a partition of unity subordinate to the cover $V_i$, the formula
$$
\si'= \sum_i \rho_i \si_i'
$$
defines  a section of $W(\om')$.\QED

Thus the set
$\Ss_\pi(M)$ of symplectic forms on
$M$, with  respect to which $\pi$ is symplectic, is open.
The aim of this paragraph is to show that a corresponding statement is true
for  Hamiltonian bundles, in other words that Hamiltonian bundles are 
{\it stable.} We begin with the following lemmas.

 \begin{lemma}\label{le:st2}  A Hamiltonian bundle $\pi:P\to B$ is stable
if and only if the restriction map
$H^2(P)\to H^2(M)$ is surjective.\end{lemma}

\proof{}   If $\pi:P\to B$ is
Hamiltonian with respect to $\om'$ then by  Proposition~\ref{prop:hamchar}
 $[\om']$ is in the image of  
$H^2(P)\to H^2(M)$.  If $\pi$ is stable, then $[\om']$ fills out a neighborhood of $[\om]$
which implies surjectivity.  Conversely, suppose that we have surjectivity.  Then 
the second condition of  Proposition~\ref{prop:hamchar} is satisfied. 
To check $(i)$ let $\ga: S^1\to B$ be a loop in $B$ and suppose that
$\ga^*(P)$ is identified symplectically with the product bundle
$S^1\times (M,\om)$.  
Let $\om_t, 0\le t\le \eps,$ be a (short) smooth path with $\om_0 =
\om$.   Then, because $P\to B$ has the structure of an $\om_t$-symplectic bundle for each
$t$, each fiber $M_b$ has a corresponding  smooth family of symplectic forms $\om_{b,t}$
of the form $g_{b,t}^* \psi_b^*(\om_t)$, where $\psi_b$ is a symplectomorphism $(M_b,
\om_b) \to (M,\om)$.  Hence, for each $t$, $\ga^*(P)$ can be symplectically
identified with $$
\bigcup_{s\in [0,1]} (\{s\}\times (M, g_{s,t}^*(\om_t))),
$$
where $g_{1,t}^*(\om_t) = \om_t$ and the $g_{s,t}$ lie in an arbitrarily small neighborhood
$U$ of the identity in $\Diff(M)$.   By Moser's homotopy argument, we can choose $U$ so
small that each $g_{1,t}$ is isotopic to the identity in the group $\Symp(M, \om_t)$.  This
proves $(i)$.\QED

\begin{lemma}\label{le:st1}\begin{itemize}
\item[(i)]
Every Hamiltonian bundle over $S^2$ is stable.
\item[(ii)]  Every symplectic bundle over a $2$-connected
base $B$ is Hamiltonian stable.
\end{itemize}
\end{lemma}
\proof{} 
$(i)$ holds because every Hamiltonian bundle over $S^2$ is c-split, in particular
the restriction map $H^2(P)\to H^2(M)$ is surjective .
$(ii)$ follows from the fact that a symplectic bundle over a
$2$-connected base is automatically Hamiltonian since
the relative homotopy groups
$\pi_{i}(\Symp(M), \Ham(M))$, $i\ge 2$, all vanish. 
(See~\cite{LM-Topology} for more details.)\QED

\begin{prop} \label{prop:stable}
Every Hamiltonian bundle is stable.
\end{prop}
\proof{}  First note that we can restrict to the case when $B$ is
simply connected.
 For the map $B \to B\,\Ham(M)$ classifying $P$ factors through a map $C\to
B\,\Ham(M)$, where $C = B/B_1$ as before, and the stability of the induced
bundle over $C$ implies that for the original bundle by Lemma~\ref{le:st2}.

Next note that by the same lemma  a
Hamiltonian bundle $P\to B$ is stable if and only if the differentials
$d_k^{0,2}: E_k^{0,2}\to E_k^{k,3-k}$ in its Leray cohomology spectral sequence vanish  
on the whole of $H^2(M)$  for $k = 2,3$. 
But it is easy to see that
we can reduce the statement for
$d_2^{0,2}$ to the case $B = S^2$.  
Thus $d_2^{0,2}= 0$ by Lemma~\ref{le:st1}$(i)$.  Similarly, we
can reduce the statement for $d_3^{0,2}$ to the case $B
= S^3$ and then use Lemma~\ref{le:st1}$(ii)$.  For more details, 
see~\cite{LM-Topology}.  \QED

It is then easy to extend the  result of \S~\ref{se:finite}
to the case when the group $\Ham(M)$ does not retract to a finite dimensional
Lie subgroup:

\begin{thm} \label{thm:stabil}  Let
$(M, \om)$ be a closed symplectic manifold, and let $\iota: X \to
\Ham(M, \om)$ be a continuous map defined on a finite 
CW-complex $X$. 
 Then, for each perturbation $\om'$ in some sufficiently small
neighbourhood $\Uu$ of $\om$ in the space of all symplectic forms on $M$,
there is a map  
$$
\iota': X \to
\Ham(M, \om')
$$
that varies continuously as the form $\om'$ varies 
in $\Uu$. 
\end{thm}

Here is the proof. Given $\iota$, construct the Hamiltonian bundle
$(M,\om) \to P \to X \times S^1$ by considering the direct 
product $(X \times [0,1]) \times (M,\om)$ and identifying 
the ends by the map
$$ 
\begin{array}{ccccc}
& \phi: & (X \times \{0\}) \times (M,\om) & \to & (X \times \{1\}) \times
(M,\om) \\ &       &   (x,0,y)                   &  \mapsto &  (x,1,\iota_x(y))
\end{array}
$$
Since $\iota$ is Hamiltonian, $\Phi$ is a Hamiltonian automorphism of
$X \times (M,\om)$ and therefore $P$ is a Hamiltonian bundle. By
Proposition~\ref{prop:stable}, the Hamiltonian structure on $P$ persists and
varies continuously as $\om$ varies in some open neighborhood $\Uu$ of
$\Ss(M)$.   For $\om' \in \Uu$, denote by $(P(\om'), \pi, \{a'\})$ the
corresponding Hamiltonian structure. Consider the restriction of $P(\om')$ to
each time $t \in S^1$. This gives a Hamiltonian bundle $(M,\om')
\hookrightarrow P_t(\om') \to  X_t = X \times \{t\}$, which for $\om' = \om$
has the trivial (i.e. product) Hamiltonian structure.  We showed in \cite{LM-Topology}
that such structures
are classified by maps $\pi_1(X_t) \to \Ga_{\om'}$ and in \cite{LMP} that the rank of
$\Ga_{\om'}$ is finite and locally constant. Thus because
the map
$\pi_1(X_t)
\to
\Ga_{\om'}$ is zero for
$\om' = \om$ and since it depends continuously on $\om'$, we conclude that
it must be zero for all $\om'$ and $t$, i.e.
the induced  Hamiltonian structure on $P_t(\om')$ is trivial for each $t$ and in
particular for $t = 0,1$.  This means that $P(\om')$ is defined as the quotient of
$(X \times [0,1]) \times (M,\om')$ by an automorphism of the trivial
Hamiltonian bundle 
$$ 
\phi_{\om'}: (X \times \{0\})\times (M,\om') \to
(X \times \{1\}) \times (M,\om'), 
$$
which is homotopic to a map $\iota': X \to
\Ham(M,\om')$.

\MS

Here is a second more direct proof. First observe that if $K$
is any compact set of $\Symp(M,\om)$, then for all $\om' \in \Ss(M)$
sufficiently close to $\om$, Moser's argument gives a canonical map
$$
\Phi_{\om,\om'}: K \to \Symp(M,\om')
$$
defined by precomposing each map $f \in K$ by 
a diffeomorphism $g_f$ of $M$ such that $g_f^*(f^*\om') = \om'$.
 Here $g$ is the diffeomorphism corresponding to
the isotopy $\om'_t = t\om' + (1-t)f^*(\om')$.  Because $K$ is compact,
one can indeed choose a small enough neighborhood $\Uu$ of $\om$ in
$\Ss(M)$ so that such a segment is nondegenerate for all $\om' \in \Uu$. 

Choosing $K$ to be the image of $\iota$, we get a map 
$\Phi_{\om,\om'} \circ \iota: B \to \Symp(M,\om')$. 
The map $\Phi_{\om,\om'} \circ \iota$
homotops into $\Ham(M,\om')$ if and only if  the elements in
$(\Phi_{\om,\om'} \circ \iota)_*(\pi_1(B))$  lie in the kernel of the
homomorphism 
$$
\Flux_{\om'} : \pi_1(\Symp_0(M,\om')) \to H^1(M,\R).
$$
In fact the flux homomorphism $\Flux_{\om'}$ is defined on $\pi_1(\Diff(M))$,
and so, since $\Phi_{\om,\om'} \circ \iota$ is homotopic to $\iota$ as maps to
$\Diff(M)$, it suffices to show that 
$\Flux_{\om'}$ vanishes on the elements of $\iota_*(\pi_1(B))$.  But  $\Flux_\om$
vanishes on $\iota_*(\pi_1(B))$ by construction, and so 
$\Flux_{\om'}$ also vanishes on these elements by
the ``stability of Hamiltonian loops"
in~\cite{LMP}.  This is just another way of expressing the stability
of Hamiltonian structures over $S^2$.  To see this, let $\phi = \iota_*(\ga)$ be the
image of a loop $\ga$ in $B$, and consider the associated bundle $P_\phi\to S^2$
constructed using $\phi$ as clutching map.  Then, for any closed form $\tau$ on
$M$, symplectic or not,  $\Flux_\tau(\phi)$ is nothing other than the value of
the Wang differential $\p_\phi$ of this bundle on the class $[\tau]$.  
The stability of $P_\phi\to S^2$ implies that $\p_\phi([\om']) = 0$, and therefore
$$
\Flux_{\om'}(\phi) = \p_\phi([\om']) = 0,
$$
as required.\MS

\subsection{From $S^2$ to more general bases, using analytic arguments} \label{ss:analytic}

\begin{prop}\label{prop:split}  Let $(M,\om)$ be a closed symplectic manifold,
and $M \hookrightarrow P \to B$ a Hamiltonian bundle 
over a CW-complex $B$. Then the
rational cohomology of $P$ splits if the base has the homotopy type of 
a symplectic manifold $W$ for which some spherical
Gromov--Witten invariant 
$n_W( pt, pt, c_1, \ldots, c_k; A)$  does not vanish,
where $k \ge 0$, $A\in H_2(W;\Z)$ and the $c_i's$ are any cycles in $W$.
\end{prop}

  Note that spaces satisfying the above condition include all products
of complex projective spaces and their blow-ups.

  A special case of this proposition was proved in \cite{LM-Topology},
and the general case will appear in \cite{L}. The proof is a
generalization  of the arguments in \cite{LMP,McD}.  The idea is to show that
 moduli spaces of $J$-holomorphic curves in ruled symplectic manifolds
$P$ behave  like fibered moduli spaces, which
implies that appropriate GW-invariants in $P$ are equal 
to the product of a GW-invariant
of the base with a GW-invariant of the fiber. 
Indeed,
suppose that $M \hookrightarrow P \to B$ is a
Hamiltonian fiber bundle over a symplectic manifold $B$ and assume that 
$B$ contains a spherical class $A$ with a 
non-zero Gromov--Witten invariant of the form 
$n_B( pt, pt, c_1, \ldots, c_k; A)$. 
Recall that this invariant counts the number of $J$\/-curves in class $A$ that
pass through two generic points and through generic representatives of
the classes $c_1,
\ldots, c_k$. 
Let $C$ be such a rational $J$\/-curve. We have explained above that the restriction
$P_C$ of $P$ to the curve $C$ is a Hamiltonian fiber bundle that c-splits. 
recall that this c-splitting was shown by taking
two $M$\/-fibers $M_0, M_{\infty}$ in $P_C$ and by finding, for each cycle $a
\in M_0$, a cycle $b = b(a) \in M_{\infty}$ such that the GW-invariant in $P_C$
$$
n_{P_C} (\iota(a), \iota(b); \si)
$$
does not vanish. (Here $\iota$ denotes the inclusion of the fiber in $P_C$,
$\si$ is some homology class of sections of $P_C \to C=S^2$ and $n_{P_C}$
counts the number of $J$\/-holomorphic curves in class $\si$ passing through
$\iota(a)$ and $\iota(b)$.)  This
implies that $\iota (a)$ cannot vanish, and therefore $P_C$ c-splits
by the Leray-Hirsch theorem. Now take an almost complex structure $J'$
on $P$ such that the projection $\pi:P \to B$ is $(J',J)$-holomorphic
and consider the invariant in $P$
$$
n_{P} (\iota'(a), \iota'(b), \pi^{-1}(c_1), \ldots, \pi^{-1}(c_k); \iota_{P_C,P}(\si))
$$
where $\pi$ is the projecton $P \to B$, $\iota'$ denotes the inclusion of the 
fiber in $P$, and $\iota_{P_C,P}$ is the inclusion of $P_C$ in $P$. It is not hard to
see, at least when the moduli spaces are well-behaved, that this last
invariant must be equal to the sum, taken over the  rational curves $C$
appearing in $n_B( pt, pt, c_1, \ldots, c_k;A)$, of the
corresponding numbers $n_{P_C} (\iota(a), \iota(b); \si)$, with signs
according to orientations in the moduli space. But because the 
Hamiltonian bundles
$P_C$ and $P_{C'}$ are isomorphic when $C$ and $C'$
are homologous in $B$, this sum is actually the product
$$
n_B(pt, pt, c_1, \ldots, c_k; A) \; \times \; n_{P_C} (\iota(a), \iota(b);
\si) 
$$
which does not vanish.

Therefore $\iota'(a)$ cannot vanish either and
the bundle $P$ c-splits.

\subsection{Iterating bundles: geometric arguments} \label{ss:Geometric-argt}

Let $M \hookrightarrow P \to
B$ be a Hamiltonian  bundle over a  simply connected base $B$ and assume that
all  Hamiltonian bundles 
over $M$ as well as over $B$ c-split. We explain in this section that 
any Hamiltonian bundle over $P$ must also be c-split.
This provides a powerful recursive argument that extends c-splitting
results to much more general bases.

We begin with some trivial observations and then discuss composites of Hamiltonian
bundles. The first lemma is true for any class of bundles with specified
 structural group.

 \begin{lemma}\label{le:pullb} Suppose that $\pi:P\to B$ is Hamiltonian and that $g:B'\to
B$ is a continuous map.  Then the induced bundle $\pi':g^*(P)\to B'$
is Hamiltonian.
\end{lemma}

Recall that any extension $\tau$ of the forms on the
fibers is called a connection form.

\begin{lemma}\label{le:sym} If $P\to B$ is a smooth Hamiltonian fiber bundle over a
symplectic base $(B,\si)$ and if $P$ is compact then  there is a connection form 
$\Om^\ka$ on $P$ that is
symplectic. \end{lemma}
\proof{}  The bundle $P$ 
carries a closed connection form $\tau$.
Since $P$ is compact, the form $\Om^\ka = \tau+\ka\pi^*(\si)$ is symplectic for large
$\ka$. \QED

Observe that the deformation type of the form $\Om^\ka$ is unique for sufficiently
large $\ka$ since given any two closed  connection forms $\tau,\tau'$   the linear isotopy
$$
t \tau + (1-t)\tau' + \ka\pi^*(\si),\quad 0\le t\le 1,
$$
consists of symplectic forms for sufficiently
large $\ka$.  However, it can happen that there is a symplectic connection
form $\tau$ such that $\tau + \ka\pi^*(\si)$ is not symplectic for small $\ka > 0$, even
though it is symplectic for large $\ka$.  (For example, suppose $P = M\times B$
and that $\tau$ is the sum $\om + \pi^*(\om_B)$ where $\om_B + \si$ is not
symplectic.)

Let us now consider the behavior of Hamiltonian bundles  under
composition. If  
$$
(M,\om)\to P
\stackrel{\pi_P}\to
X,\quad\mbox{and}\quad (F,\si)\to X\stackrel{\pi_X}{\to} B
$$
 are Hamiltonian
fiber bundles,  then the restriction 
$$
\pi_P:\quad W = \pi_P^{-1}(F)\longrightarrow F
$$
is a Hamiltonian fiber bundle.  Since $F$ is a manifold, we can assume
without loss of generality that $W\to F$ is smooth.  
Moreover, the manifold $W$
carries a symplectic connection form $\Om_W^\ka$, and
 it is natural to ask when the composite
map $\pi:P\to B$ with fiber $(W, \Om_W^\ka)$  is itself Hamiltonian.   

\begin{lemma}\label{le:comp}  Suppose
 that $B$ is a simply connected CW-complex and that $P$ is compact.
 Then  $\pi = \pi_X\circ \pi_P: P\to B$ is a 
Hamiltonian fiber bundle with fiber $(W, \Om_W^\ka)$, where
$ \Om_W^\ka = \tau_W  + \ka\pi_P^*(\si)$,
$\tau_W$ is any symplectic connection form  on $W$, and $\ka$ is sufficiently large.
 \end{lemma} \proof{} 
We may assume that the base $B$ as well as the bundles
are smooth. Let $\tau_P$ (resp.
$\tau_X$) be a closed connection form with respect to the bundle $\pi_P$, (resp. $ \pi_X$),
and let $\tau_W$ be its restriction to $W$.
Then $\Om_W^\ka$ is the restriction to $W$ of the closed form
$$
\Om_P^\ka = \tau_P + \ka\pi_P^*(\tau_X).
$$
By increasing $\ka$ if necessary we can ensure that $\Om_P^\ka$ restricts to a symplectic
form on every fiber of $\pi$ not just on the the chosen fiber $W$.
 This shows firstly that
$\pi:P\to B$ is symplectic, because there is a well defined 
 symplectic form on each of its fibers, and secondly that it is Hamiltonian with respect to
this form $\Om_W^\ka$ on the fiber $W$.   Hence 
Lemma~\ref{le:st2} implies that  $H^2(P)$ surjects onto $H^2(W)$. 

Now suppose that  $\tau_W$ is {\it any}
closed connection form on $\pi_P: W\to F$.  
Because the restriction map $H^2(P) \to H^2(W)$ is surjective,
the cohomology class $[\tau_W]$ is the restriction of a class on $P$
and so, by  Thurston's construction, the form $\tau_W$ can be extended to a closed connection
form
$\tau_P$  for the bundle $\pi_P$. 
Therefore the previous argument applies in this case too. \QED

Now let us consider the general situation, when $\pi_1(B) \ne 0$.  
The proof of the lemma above applies to show that the composite bundle $\pi:P\to B$ 
is symplectic with respect to suitable $\Om_W^\ka$ and that it has a symplectic connection
form.  However,
 even though $\pi_X:X\to B$ is symplectically trivial over the $1$-skeleton $B_1$ the
same may not be true of the composite map $\pi:P\to B.$  Moreover, in general it is not
clear whether triviality with respect to one form $\Om_W^\ka$ implies that for another.
Therefore, we may conclude the following:

\begin{prop}\label{prop:comp}
If
$(M,\om)\to P
\stackrel{\pi_P}\to
X,$ and $ (F,\si)\to X\stackrel{\pi_X}{\to} B$ are Hamiltonian fiber bundles and $P$ is
compact, then the composite
 $\pi =\pi_X\circ\pi_P:P\to B$ is a symplectic fiber bundle with respect to any form
$\Om_W^\ka$ on its fiber $W = \pi^{-1}(pt)$, provided that $\ka$ is sufficiently large. 
Moreover if $\pi$ is symplectically trivial over the $1$-skeleton of $B$ with respect to
 $\Om_W^\ka$ then $\pi$ is Hamiltonian.
\end{prop}

In practice, we will apply these results in cases 
where $\pi_1(B) = 0$.  We will not specify
the precise form on $W$, assuming that 
it is $\Om_W^\ka$ for a suitable $\ka$.

\begin{lemma}\label{le:funct2}
If $(M,\om)\stackrel{\pi}\to  P\to B$ is a 
compact Hamiltonian  bundle over a simply connected CW-complex
$B$ and if every Hamiltonian fiber bundle over 
$M$ and $B$ is c-split, then every Hamiltonian
bundle over $P$ is c-split. 
\end{lemma}
\proof{}  Let $\pi_E: E\to P$ be a Hamiltonian bundle with fiber $F$ and let
$$
F\to W\to M
$$
be its restriction over $M$.    Then by assumption the latter bundle
c-splits so that $H_*(F)$ injects into $H_*(W)$.  Lemma~\ref{le:comp} 
implies that
the composite bundle $E\to B$ is
Hamiltonian with fiber $W$ and therefore also c-splits.   Hence $H_*(W)$
injects into $H_*(E)$. Thus $H_*(F)$ injects into $H_*(E)$, as required.\QED

\subsection{Topological arguments} \label{ss:Topological-argt}

We now put together the results and methods of the last 
subsections about c-splitting.  For more details see \cite{LM-Topology}.

\begin{lemma} \label{le:sur} If $\Si$ is a
closed orientable surface then any Hamiltonian bundle over 
$S^2 \times \ldots \times S^2 \times \Si$
 is c-split. 
\end{lemma}  
\proof{} Consider any degree one map $f$ from $\Si\to S^2$. 
Because $\Ham(M, \om)$ is
connected, $B\,\Ham(M, \om)$ is simply connected, and therefore any
homotopy class of maps from $\Si\to B\,\Ham(M, \om)$ factors through
$f$.  Thus any Hamiltonian bundle over $\Si$ is the pullback by $f$ of a Hamiltonian
bundle over $S^2$. Because such bundles  c-split
 over $S^2$, the same is true over $\Si$ by Lemma~\ref{le:funct}$(i)$.

The statement for $S^2 \times \ldots \times S^2 \times \Si$
 is now a direct consequence of iterative 
applications of Lemma~\ref{le:funct2} applied to the trivial bundles 
$
 S^2 \times \ldots \times S^2 \times \Si \to S^2.
$
\QED

\begin{cor} \label{cor:odd} Any Hamiltonian bundle over
$S^2 \times \ldots \times S^2 \times S^1$ is c-split.
\end{cor}

\proof{} Consider the maps $S^1 \to T^2 \to S^1$ given by inclusion on
the first factor and projection onto the first factor. Their composition
is the identity. Extend them to maps 
$$
S^2 \times \ldots \times S^2 \times S^1 \to 
    S^2 \times \ldots \times S^2 \times T^2  \to
             S^2 \times \ldots \times S^2 \times S^1 .
$$
by multiplying with the identity on the $S^2$ factors. Then a
 Hamiltonian bundle $P$ on $S^2 \times \ldots \times S^2 \times S^1$
pulls-back to a c-split bundle $P'$ on 
$S^2 \times \ldots \times S^2 \times T^2$ by 
Lemma~\ref{le:sur}. By naturality, its pull-back $P''$
to $S^2 \times \ldots \times S^2 \times S^1$ is c-split. But $P'' = P$.
\QED

\begin{prop}\label{prop:sph} For each $k\ge 1$, every 
Hamiltonian bundle over $S^k$  c-splits.
\end{prop}
\proof{}  By Lemma~\ref{le:sur} and 
Corollary~\ref{cor:odd} there is for each
$k$ a  $k$-dimensional closed manifold 
$X$ such that every Hamiltonian bundle 
over $X$ c-splits.  Given any Hamiltonian
bundle $P\to S^k$ consider its pullback to 
$X$ by a map $f:X\to S^k$ of degree $1$.
Since the pullback c-splits, the original bundle does too by
Lemma~\ref{le:funct}$(ii)$.  \QED

 By the Wang exact sequence, this implies that the action of 
the homology groups of $\Ham(M)$ on $H_*(M)$ is always trivial.

Here are some other examples of situations in which Hamiltonian
bundles are c-split.

\begin{lemma} \label{le:products} Every Hamiltonian 
bundle over $\CP^{n_1} \times \ldots \times \CP^{n_k}$ c-splits.
\end{lemma}
\proof{}
This is an obvious application of  Lemma~\ref{le:funct2}.
\QED

\begin{lemma} \label{le:three} Every Hamiltonian 
bundle over a compact CW-complex of dimension $\le 3$ c-splits.
\end{lemma}
\proof{}
This is because one can first assume that $B$ is simply-connected
and  then construct a homology surjection
$B' \to B$ where $B'$ is a wedge of $2$ and $3$\/-spheres.
\QED

 \begin{prop}\label{le:spheres} Every Hamiltonian bundle over
a product of spheres c-splits, provided that
there are no more than $3$ copies of $S^1$. \end{prop}
\proof{}  By hypothesis $B = \prod_{i\in I} S^{2m_i}\times \prod_{j\in
J}S^{2n_i + 1}\times T^k$,  where $n_i > 0$ and $0\le k\le 3$.  Set
 $$
B' =  \prod_{i\in I} \CP^{m_i} \times \prod_{j\in J} {\CP}^{n_i} \times T^{|J|}
\times T^\ell,
$$
where $\ell = k$ if $k + |J|$ is even and $= k+1$ otherwise.
Since $\CP^{n_i}\times S^1$ maps onto $S^{2n_i+1}$ by a map of degree $1$,
there is a homology surjection $B'\to B$ that maps the factor $T^\ell$ to
$T^k$. By the surjection lemma, it suffices to show  that  the pullback bundle
$P'\to B'$ is c-split. 

Consider the fibration 
$$
T^{|J|}\times T^\ell \to B' \to 
\prod_{i\in I} \CP^{m_i} \times \prod_{j\in J} {\CP}^{n_i}.
$$
Since $|J| + \ell$ is even, we can think of this as a Hamiltonian bundle.
Moreover, by construction, the restriction of the bundle $P'\to B'$ to
$T^{|J|}\times T^\ell$ is the pullback of a bundle over $T^k$, since the map
$T^{|J|}\to B$ is nullhomotopic. (Note that each $S^1$ factor in  $T^{|J|}$ goes
into a different sphere in $B$.)   Because $k\le 3$, the bundle over $T^k$ c-splits.
Hence we can apply the argument in Lemma~\ref{le:funct2} to conclude that 
$P'\to B'$ c-splits.\QED

\begin{cor} Let $B$ be a simply connected Lie group, or more generally
any  H-space whose rational fundamental group
has rank less than $4$ and whose homotopy groups are finitely 
generated in each dimension. Then c-splitting holds for all Hamiltonian
bundles over  $B$.
\end{cor}
\proof{}  Let $B$ be such a $H$-space. By the theory of minimal
models (see \cite{DMV} for instance) which applies in this case
because the fundamental group of $B$ acts trivially on all
higher homotopy groups, the rational cohomology
of $B$ is generated as a $\Q$-vector space by cup-products
of elements that pair non-trivially with spheres, ie
each $a \in H^*(B,\Q)$ can be written as a cup product
$\cup_i a_i$'s where there is for each $i$ a 
spherical class $\al_i$ in rational homology with
$a_i(\al_i) \neq 0$. If we denote by the same
symbol $\al_i: S^{n_i} \to B$ a map that realises a
non-zero multiple of the class $\al_i$, then the obvious map
$\vee_i \al_i: \vee_i S^{n_i} \to B$ extends to a map $\phi_a$ defined on the
product of these spheres that pulls back the element $a$ to
a generator of the top rational cohomology group. If there were a Hamiltonian
bundle $P$ over $B$ that did not c-split, there would be an element of lowest
degree
$a \in H^*(B;\Q)$ with non-zero differential in the spectral sequence of $P$ 
and therefore the differential of the corresponding top element of 
$H^*(\Pi_i S^{n_i})$ in the spectral sequence of the
pull back bundle $\phi_a^*(P)$ would not vanish either. But this
contradicts the c-splitting established in the previous proposition.
\footnote{We are grateful to Jaroslaw Kedra who pointed out a variant of this
argument to us.}
\QED

In particular,  Hamiltonian fibrations c-split over the loop space 
$\Om X$ of any simply connected CW-complex $X$ with $\pi_{2}X$ of small 
enough rank.  It is not at all clear how to go from this fact to c-splitting 
over $X$.  The paper~\cite{LM-Topology} contains an extensive 
discussion of what can be proved when $X$ has dimension $4$.

\begin{lemma}\label{le:bott}  Every
Hamiltonian bundle over a coadjoint orbit c-splits.
\end{lemma}
\proof{}  This follows immediately
from the remarks in Grossberg--Karshon \cite{GK}\S3 about Bott towers.  A
Bott tower is an iterated bundle $M_k\to M_{k-1} \to\dots \to M_1 = S^2$ of
K\"ahler manifolds where each map $M_{i+1}\to M_i$ is a bundle with fiber
$S^2$.  They show that any coadjoint orbit $X$ can be blown up to a
manifold that is diffeomorphic to a Bott tower $M_k$.  Moreover the
blowdown map $M_k\to X$ induces a surjection on
rational homology.   Every Hamiltonian bundle over $M_k$ c-splits by repeated
applications of Lemma~\ref{le:funct2}. Hence the result follows from the
surjection lemma.\QED

\section{Applications to ruled symplectic manifolds}

\begin{thm}\label{thm:obstruction} $\,$ {\rm \bf Obstructions to the existence of
ruled symplectic structures.} Let $M$ be a closed manifold and $P$ a
smooth fiber bundle with fiber $M$ over a simply connected manifold $B$ and
assume that $B$ is either a compact CW-complex of dimension less than $4$
or is a product of complex projective
spaces, of spheres and of coadjoint orbits of arbitrary dimensions. 
Denote by $\iota$ the inclusion of
the fiber in
$P$. Then the non-vanishing of the kernel of
$$
\iota_{\ast}: H_{\ast}(M) \to  H_{\ast}(P)
$$
is an obstruction to the existence of a ruled symplectic structure on $P$.
\end{thm}

By the Leray-Hirsch theorem, the vanishing of the kernel in the theorem above
amounts to the cohomological splitting $H_*(P) = H_*(B) \otimes H_*(M)$. 
Thus this last result may be stated as follows: under the given conditions
on $B$, a ruled structure exists on $P$ only if $P$ splits cohomologically.
This imposes strong topological constraints on the construction of
ruled symplectic manifolds by twisted products of two given ones.

Theorem~\ref{thm:obstruction} is an immediate
corollary of our results about c-splitting and of the characterization
of Hamiltonian fiber bundles over simply connected bases
in terms of the existence of a closed
extension to the total space of the symplectic forms on the fibers. 
 This
characterization also implies the following version of Hamiltonian stability.

\begin{thm} \label{thm:stability-bundles} $\,$ {\rm \bf Stability
of ruled  symplectic structures.}
Let $M \hookrightarrow P \stackrel{\pi}{\to} B$ be a smooth compact
fiber bundle over a simply connected
manifold $B$. Suppose that $P$ admits a ruled symplectic
structure $\Om$, that restricts to $\om$ on the $M$\/-fiber. Then
the ruled symplectic structure on $P$ persists under small deformations of
$\om$, i.e. there is a neighborhood $\Uu$ of $\om$ in the space of all
symplectic forms on $M$ such that each $\om' \in \Uu$ extends to a ruled
symplectic structure $\Om'$ on $P$, which varies continuously as $\om'$ varies
in $\Uu$.
\end{thm}

Observe that the above theorem remains
true for arbitrary bases $B$ provided that $P\to B$ is 
symplectically trivial over
the $1$-skeleton of $B$.

\section{Concluding remarks}

  It is still unclear whether every Hamiltonian fiber
bundle over any compact CW-complex c-splits. One of the simplest unknown
cases is a Hamiltonian bundle $(M,\Om) \to P \to B$ with
base the $4$\/-torus and with fiber a symplectic $4$\/-manifold
that does not satisfy specific properties like 
the hard Lefschetz property. The problem here
is that, if one tries to apply Lemma~\ref{le:funct2}
to $(M,\om) \to P \stackrel{\rho}{\to} B$ with $B = T^4$ given itself as
a bundle  $T^2 \to T^4 
\stackrel{\pi}{\to} T^2$, then nothing guarantees 
that the composite fibration $W
\to P \stackrel{\pi\circ \rho}\to T^2$ is trivial over the $1$-skeleton of the base.
In fact, the structural group of the  composite fibration may well be
a disconnected 
subgroup of the symplectomorphism group of the fiber $W = 
(\pi \circ \rho)^{-1}(pt)$.  Note, however,  that because 
all Hamiltonian fibrations over $T^3$ c-split, 
we do know that the elements of this subgroup act trivially on the 
cohomology of $W$.  This raises the interesting question of whether one can
extend our results on c-splitting for Hamiltonian bundles to certain
disconnected extensions of the Hamiltonian group.
We have no techniques at present to deal with this question, since bundles over
$T^2$ need not admit any $J$-holomorphic sections.

\end{document}